\documentclass[12pt]{article}
\usepackage{amsmath, amsthm, amssymb, fullpage}

\numberwithin{equation}{section}
\newtheorem{theorem}[equation]{Theorem}
\newtheorem{lemma}[equation]{Lemma}
\newtheorem{cor}[equation]{Corollary}
\newtheorem{prop}[equation]{Proposition}
\theoremstyle{definition}
\newtheorem{defn}[equation]{Definition}
\newtheorem{example}[equation]{Example}
\newtheorem{remark}[equation]{Remark}
\newtheorem{convention}[equation]{Convention}

\def\be{\mathbf{e}}
\def\bv{\mathbf{v}}
\def\bw{\mathbf{w}}
\def\bx{\mathbf{x}}

\newcommand{\QQ}{\mathbb{Q}}
\newcommand{\RR}{\mathbb{R}}
\newcommand{\ZZ}{\mathbb{Z}}

\newcommand{\gotho}{\mathfrak{o}}

\DeclareMathOperator{\charac}{char}
\DeclareMathOperator{\dist}{dist}
\DeclareMathOperator{\ndist}{ndist}

\begin{document}

\title{Basis discrepancies for extensions of valued fields}
\author{Kiran S. Kedlaya}
\date{April 9, 2006}

\maketitle

\begin{abstract}
Let $F$ be a field complete for a real valuation.
It is a standard result in valuation theory that a finite extension of $F$
 admits a valuation basis if and only if it is without defect.
We show that even otherwise, one can construct bases in which the discrepancy
between measuring valuation an element versus on the components
in its basis decomposition can be made arbitrarily small.
The key step is to verify this for extensions of degree equal to the 
characteristic by a direct calculation.
\end{abstract}

\section{Introduction}

When working with valuations of fields, it is often useful to be able 
to calculate valuations on an extension field in terms of valuations on the
base field. In particular, if $E/F$ is a finite extension of fields,
$F$ is equipped with a real valuation $v$ that extends uniquely to $E$
(e.g., $F$ is henselian), and $\be_1, \dots, \be_n$ is a basis for $E$ as
a vector space over $F$, one can ask whether the inequality
\[
v(a_1 \be_1 + \cdots + a_n \be_n) \geq \min_i \{v(a_i \be_i)\}
\]
is always an equality; this would be the valuation-theoretic analogue of
the $\be_i$ forming an orthogonal basis.

Sadly, it is not always possible to find such a basis for a given field
extension, and the culprit is the usual one: it is the presence of
``defect'', i.e., the failure of the equality between total degree and
the product of ramification and inertial degrees. However, the discrepancy
\[
v(a_1 \be_1 + \cdots + a_n \be_n) - \min_i \{v(a_i \be_i)\}
\]
turns out to be bounded for any particular basis, and one can ask whether
one can at least make it arbitarily small by making good choices of basis.
We show (Theorem~\ref{T:almost val basis}) that this is indeed the case.

The key calculation for the proof of Theorem~\ref{T:almost val basis}
is the computation of discrepancy for a monomial
basis in an extension of degree equal to the characteristic
(Proposition~\ref{P:degree p}). This result may be useful in its own right;
we suspect it can be generalized to higher degree at least in the case of a basis
of a purely inseparable extension generated by a $p$-basis, though we are
only able to partially verify this (Section~\ref{sec:p-bases}).

\subsection*{Acknowledgments}
Thanks to Franz-Viktor Kuhlmann for helpful discussions, and for making
his extremely useful book in progress \cite{kuhlmann-book} available online.
The author was supported by NSF grant DMS-0400747.

\section{Valuations}

We first set some notations.

\begin{defn}
Let $G$ be an abelian group.
A \emph{(real) valuation} on $G$ is a function $v: G 
\to \RR \cup \{\infty\}$ with the following properties.
\begin{enumerate}
\item[(i)] For $x \in G$, $x=0$ if and only if $v(x) = +\infty$.
\item[(ii)] For $x,y \in G$, $v(x+y) \geq \min\{v(x), v(y)\}$.
\end{enumerate}
We say the valuations $v,w: G \to \RR \cup\{\infty\}$ are \emph{metrically
equivalent} if there exist $c, d \in \RR$ such that for all $x \in G$,
\[
v(x) + c \leq w(x) \leq v(x) + d;
\]
this is clearly an equivalence relation on valuations, and two metrically
equivalent valuations induce the same ultrametric topology on $G$.
By a \emph{valued (abelian) group}, we mean an abelian 
group $G$ equipped with a real valuation $v$.
\end{defn}

\begin{defn}
Let $F$ be a field. A \emph{(real) Krull valuation} on $F$ is a valuation $v$
on the additive group of $F$, which satisfies the following
additional property.
\begin{enumerate}
\item[(iii)]
For $x,y \in F$, $v(xy) = v(x) + v(y)$.
\end{enumerate}
By a \emph{valued field}, we mean a field $F$ equipped with a real 
Krull valuation 
$v$. For $(F,v)$ a valued field, let $\gotho_F$ and $\kappa_F$ denote the
valuation ring and residue field, respectively,
and let $\hat{F}$ be the completion of $F$, which is a field to which $v$
extends uniquely.
\end{defn}

We defer to \cite{vaquie} for a more detailed discussion of Krull valuations
and their properties.

\begin{defn}
Let $(F, v)$ be a valued field, and let $V$ be an $F$-vector space. 
An \emph{$F$-valuation} on $V$ is a valuation $f: V \to \RR \cup \{\infty\}$
on the underlying additive group of $V$, which satisfies the following
additional property.
\begin{enumerate}
\item[(iii)]
For $x \in F$ and $\bv \in V$, $f(x\bv) = v(x) + f(\bv)$.
\end{enumerate}
By a \emph{valued $F$-vector space}, we will
mean an $F$-vector space $V$ equipped with a valuation $f$.
\end{defn}

\begin{example}
Let $(F,v)$ be a valued field, and let $E$ be a finite extension of $F$.
Then the set of extensions of $v$ to a Krull valuation on $E$ 
is nonempty and finite; in particular,
it consists of a single element in case $F$ is henselian, or if
$E/F$ is purely inseparable.
\end{example}

\section{Discrepancies}

\begin{defn}
Let $F$ be a valued field, let $V$ be a finite dimensional
valued $F$-vector space,
and let $T = \{\be_1, \dots, \be_n\}$ be
an $F$-basis of $V$. 
For $\bv \in V$, write $\bv = c_1 \be_1 + \cdots + c_n \be_n$,
and  define the \emph{discrepancy} of
$T$ with respect to $\bv$ to be
\[
d(T; \bv) = f(\bv) - \min_i \{ f(c_i \be_i)\};
\]
note that this quantity is always nonnegative.
Define the \emph{discrepancy} of $T$ to be the initial segment
\[
d(T) = \sup_{\bv \neq 0} \{ d(T; \bv)\};
\]
note that $d(T)$ could be infinite in general, but also note
Corollary~\ref{C:discrepancy} below.
If $d(T) = 0$, we say that $T$ is a \emph{valuation basis} of $V$ over $F$.
\end{defn}

\begin{remark}
The discrepancy is additive in the following sense.
Let $F$ be a valued field, and suppose $E$ is a finite extension of $F$
to which $v$ admits a unique extension; use this extension to view
$E$ as a valued $F$-vector space. Let $V$ be a finite dimensional valued
$E$-vector space; then $V$ is also a valued $F$-vector space.
If $T$ is a basis of $E$ over $F$, and $S$
is a basis of $V$ over $E$, then
\[
ST = \{st: s \in S, t \in T\}
\]
is a basis of $V$ over $F$, and we have $d(ST) = d(S) + d(T)$.
\end{remark}

\begin{lemma} \label{L:equiv norms}
Let $F$ be a complete valued field, and let $V$ be a finite dimensional
$F$-vector space.
Then any two $F$-valuations on $V$ are metrically equivalent.
\end{lemma}
\begin{proof}
See \cite[Lemma~2]{roquette}.
\end{proof}
\begin{cor} \label{C:discrepancy}
Let $F$ be a complete valued field, and let $V$ be a finite dimensional
$F$-vector space. Then for any $F$-basis $T$ of $V$, the discrepancy
$d(T)$ is finite.
\end{cor}
\begin{cor}
Let $(F,v)$ be a valued field, and suppose that $E$ is a finite
extension of $F$ linearly disjoint from $\hat{F}$;
use this extension to view $E$ as a valued $F$-vector space.
Then for any $F$-basis $T$ of $E$, the discrepancy
$d(T)$ is finite.
\end{cor}
\begin{proof}
By hypothesis $E \otimes_F \hat{F}$
is a field, so the map from it to $\hat{E}$ is injective. Since these
are $\hat{F}$-vector spaces of the same dimension, they are actually isomorphic.
We may thus apply Corollary~\ref{C:discrepancy} to obtain the finiteness
of the discrepancy of $T$ as a basis of $\hat{E}$ over $\hat{F}$;
the discrepancy of $T$ as a basis of $E$ over $F$ cannot be any larger.
\end{proof}

\section{Valuation bases and defect}

\begin{defn}
Let $F$ be a valued field, and let $V \subseteq W$ be an inclusion of
finite dimensional valued $F$-vector spaces. Let $\be_1, \dots, \be_n$ be
a basis of $W$ over $V$, i.e., a sequence of elements of $W$ which
form a basis in $W/V$. We say $\be_1, \dots, \be_n$ is a \emph{valuation
basis} of $W$ over $V$ if for any $\bw \in W$, if we write
$\bw = \bv + c_1 \be_1 + \cdots + c_n \be_n$ with $\bv \in V$
and $c_i \in F$, we have
\[
v(\bw) = \min\{v(\bv), v(c_1 \be_1), \dots, v(c_n \be_n)\}.
\]
Note that if $V = \{0\}$, then a valuation basis of $W$ over $V$ is the same
as a valuation basis of $W$ in our previous sense.
\end{defn}

\begin{lemma} \label{L:val bases}
Let $F$ be a valued field, and let $V \subseteq W \subseteq X$ be inclusions of
finite dimensional valued $F$-vector spaces. Suppose that $X$ admits
a valuation basis over $V$. Then $X$ also admits a valuation basis over $W$.
\end{lemma}
\begin{proof}
By induction on $\dim(W/V)$, it suffices to consider the case where
$\dim(W/V) = 1$. Choose $\bw \in W \setminus V$, 
let $\be_1, \dots, \be_n$ be a valuation
basis of $X$ over $V$, and write $\bw = \bv + c_1 \be_1 + \cdots + c_n \be_n$
for some $\bv \in V$ and some $c_1, \dots, c_n \in F$.
Choose $j \in \{1, \dots, n\}$ to minimize $v(c_j \be_j)$; then
$v(\bw - \bv) = v(c_j \be_j)$ because the $\be_i$ form a valuation basis. 

We claim that
omitting $\be_j$ yields a valuation basis for $X$ over $W$. Namely,
for $\bx \in X$, we can write $\bx = \bv' + d (\bw-\bv) 
+ \sum_{i \neq j} d_i \be_i$
for some $\bv'\in V$ and some $d,d_i \in F$. We also have
\[
\bx = \bv' + d c_j \be_j + \sum_{i \neq j} (d_i + d c_i) \be_i
\]
which yields
\[
v(\bx) = \min\{v(\bv'), v(dc_j \be_j), \min_{i \neq j} \{ v((d_i + d c_i) \be_i)
\}.
\]
If this last minimum is achieved by $v(\bv')$ or $v(dc_j \be_j) = 
v(d(\bw - \bv))$, we are done. Otherwise, for any $i$ achieving
$\min_{i \neq j} \{v((d_i + d c_i)\be_i)\}$, we have
\[
v(d c_i \be_i) \geq v(d c_j \be_j)
> v((d_i + d c_i) \be_i)
\]
and so $v((d_i + d c_i) \be_i) = v(d_i \be_i)$, and we are done again.
\end{proof}

\begin{defn}
Let $(F,v)$ be a henselian valued field, and let $E$ be a finite
extension of $F$, to which $v$ necessarily extends uniquely. Define the
\emph{ramification degree} of the extension $E/F$ as the group index
$e_{E/F} = [v(E^*):v(F^*)]$, and the \emph{inertia degree} as the field degree
$f_{E/F} = [\kappa_E:\kappa_F]$.
\end{defn}

The key relationship between degree, ramification degree, and 
inertia degree is the following lemma of Ostrowski.
\begin{convention}
For $F$ a field, define the \emph{characteristic exponent} of $F$ to be 1
if $\charac(F) = 0$ and $\charac(F)$ otherwise.
\end{convention}

\begin{prop}[Ostrowski] \label{P:Ostrowski}
Let $(F,v)$ be a henselian valued field of characteristic exponent $p$.
Let $E$ be a finite extension of $F$.
Then
there exists a nonnegative integer $\delta_{E/F}$ such that
\[
[E:F] = p^{\delta_{E/F}} e_{E/F} f_{E/F}.
\]
\end{prop}
\begin{proof}
See
\cite[Th\'eor\`eme~2, p.\ 236]{ribenboim}.
\end{proof}

\begin{defn}
With notation as in Proposition~\ref{P:Ostrowski}, we refer to
$p^{\delta_{E/F}}$ as the \emph{defect} of the extension $E/F$.
We say $E/F$ is \emph{defectless} (or \emph{without defect})
if the defect is equal to 1.
\end{defn}

\begin{example}
For the benefit of the reader unfamiliar with the concept of defect, it
is worth recalling a simple example where it arises, from
\cite[\S 3.3]{kuhlmann-value}.
Let $k$ be an algebraically
closed field of characteristic $p>0$. Let $e_1, e_2, \dots$ be a sequence of
positive integers with $e_{i+1} \geq e_i + i$ for all $i$. Define
\[
y = \sum_{i=1}^\infty x^{-p^{-e_i}}
\]
in the ring $k((x^{\QQ}))$ 
of Hahn-Mal'cev-Neumann series over $k$ in the variable $x$.
Put $F = k(x,y) \subset k((x^{\QQ}))$ and let $v$ be the restriction to $F$
of the $x$-adic valuation on $k((x^{\QQ}))$. One then verifies that the
group $v(F^*)$ is $p$-divisible, so the ramification index of any finite
extension of $F$ is coprime to $p$. In particular, the extension
\[
E = F[z]/(z^p - z - x^{-1})
\]
has degree $p$, so $e_{E/F} = 1$; also $f_{E/F} = 1$ because $k$ is already
algebraically closed, so $\delta_{E/F} = [E:F] = p$.
\end{example}

Although the following result seems to be folklore, we will make it explicit
here.
\begin{prop} \label{P:defect basis}
Let $(F,v)$ be a henselian valued field, and let $E$ be a finite
extension of $F$.
Then $E/F$ is defectless
if and only if $E$, viewed as a valued $F$-vector space, admits a valuation
basis.
\end{prop}
\begin{proof}
Suppose $E/F$ is defectless. Choose elements $x_1, \dots, x_m \in \gotho_E$
lifting a basis for $\kappa_E$ over $\kappa_F$, and choose elements
$y_1, \dots, y_n \in E$ whose images under $v$ represent the cosets
of $v(F^*)$ in $v(E^*)$. Put 
\[
T = \{x_i y_j: i = 1, \dots, m; \, j = 1, \dots, n\};
\]
then $T$ is a valuation basis of its $F$-span within $E$. But that span
has dimension $e_{E/F} f_{E/F} = [E:F]$ since $E/F$ is defectless, so $E/F$
admits a valuation basis.

Conversely, suppose $E/F$ admits a valuation basis. Construct $T$ as above,
and let $V$ be its $F$-span. Suppose by way of contradiction that $E \neq V$.
By Lemma~\ref{L:val bases}, $E$ admits a valuation
basis $\be_1, \dots, \be_n$ over $V$ as vector spaces over $F$.
Pick any $\bw \in E \setminus V$; then
$\sup_{\bv \in V} \{v(\bw - \bv)\}$ is achieved, because we can write
$\bw = \bv + c_1 \be_1 + \cdots + c_n \be_n$ for some $\bv \in V$
and some $c_i \in F$, and this $\bv$ works. However, by the way we constructed
$T$, for any $\bv \in V$, we can construct $\bv' \in V$ such that
$v(\bw - \bv - \bv') > v(\bw - \bv)$, so the supremum cannot be achieved.
This contradiction yields the equality $E = V$, so $E/F$ is defectless as 
desired.
\end{proof}
\begin{cor} \label{C:tame val basis}
Let $(F,v)$ be a henselian valued field, and let $E$ be a finite
extension of $F$.
If the degree $[E:F]$
is coprime to the characteristic exponent of $F$, then $E/F$ admits
a valuation basis.
\end{cor}

\section{Distance and discrepancy}

\begin{defn}
Let $G$ be a valued group. For $x \in G$ and $S \subseteq G$, we define
the \emph{distance from $x$ to $S$} as 
\[
\dist(x, S) = \sup\{v(x-y): y \in S\}.
\]
If $S$ is a subgroup, we define the \emph{normalized distance}
from $x$ to $G$ as
\[
\ndist(x, S) = \dist(x,S) - v(x) = \sup\{v(x-y) - v(x): y \in S\}.
\]
Note that we can also write
\[
\ndist(x,S) = \sup\{v(x-y) - \min\{v(x),v(y)\}: y \in S\},
\]
since we can ignore the contributions from those $y$ with $v(y) < v(x)$
(they give nonpositive quantities in the first sup and zero in the second sup,
whereas $\ndist(x,S) \geq 0$ because we can put $y=0$).
\end{defn}

\begin{prop} \label{P:degree p}
Let $(F,v)$ be a henselian valued field of characteristic $p>0$. Let $E$ be an extension
(separable or not)
of $F$ of degree $p$, and suppose
$x \in E$ generates $E$ over $F$. Put $T = \{1,x,\dots,x^{p-1}\}$, viewed
as a basis of $E$ over $F$. Then
\[
d(T) = (p-1)\ndist(x,F).
\]
\end{prop}
\begin{proof}
We first show that $d(T) \geq (p-1) \ndist(x,F)$. 
For any $y \in F$, we have
\begin{align*}
d(T; (x-y)^{p-1}) &= d(T; \sum_{i=0}^{p-1} \binom{p-1}{i} x^{p-i} (-y)^i) \\
&= (p-1)v(y-x) - \min_{0\leq i \leq p-1} \{(p-1-i) v(x) + i v(y)\} \\
&= (p-1)v(y-x) - (p-1) \min\{v(x), v(y)\}.
\end{align*}
The supremum of the last expression is
$(p-1) \ndist(x,F)$, while the supremum of
$d(T; (x-y)^{p-1})$ is bounded above by $d(T)$, whence the desired inequality.

We next show that $d(T) \leq (p-1) \ndist(x,F)$. For a nonzero polynomial
$P(t) = \sum a_i t^i$ over $F$, define
\[
w(P) = \min_i \{v(a_i) + i v(x)\}
\]
and note that this gives a valuation on $F[t]$.

Now fix a polynomial $P(t) \in F[t]$ of degree $d \leq p-1$.
Let $F'$ be a splitting field of $F$;
then $[F':F]$ is coprime to $p$, so $F'$ and $E$ are linearly
disjoint. Let $E'$ be the compositum of $F'$ and $E$
over $F$, so that $E'$ is an extension of $F'$ of degree $p$ generated by $x$.
By Corollary~\ref{C:tame val basis}, $F'$ admits a valuation basis as an
$F$-vector space, which we can rescale to force 1 into it,
and the same basis serves as a valuation basis for $E'$ as an $E$-vector
space.
This implies that for any $y \in E$,
\[
\dist(y, F') = \dist(y, F).
\]
Over $F'$, we can factor $P(t)$ as a product $a_d \prod_{i=1}^d (t-r_i)$, and compute
\begin{align*}
d(T; P(x)) &= v(P(x)) - w(P) \\
&= \sum_{i=1}^d (v(x-r_i) - \min\{v(x), v(r_i)\}) \\
&\leq d \ndist(x, F') \\
&\leq (p-1) \ndist(x, F),
\end{align*}
yielding the desired inequality and completing the proof.
\end{proof}

{}From this simple calculation we obtain a striking conclusion.

\begin{theorem} \label{T:almost val basis}
Let $(F,v)$ be a henselian valued field, and let $E$ be a finite extension of 
$F$.
Then for any $c \in \Gamma_{> 0}$, there
exists a basis of $E$ over $F$ with discrepancy less than $c$.
\end{theorem}
\begin{proof}
We proceed by induction on $[E:F]$; by virtue of the additivity of discrepancy,
if at any point we can insert an intermediate field between $E$ and $F$,
we may reduce to considering the two intermediate extensions.

First of all, we can always insert some $E'$ between $E$ and $F$ so that
$E'/F$ is separable and $E/E'$ is purely inseparable. Since a purely inseparable
extension can be written as a tower of extensions of monogenic extensions
of degree $p$, we can apply Proposition~\ref{P:degree p} successively to
treat that case. We may thus assume hereafter that $E/F$ is separable.

Since $E/F$ is a finite separable extension of henselian valued fields,
there exists a tower of field extensions
\[
F \subseteq U \subseteq T \subseteq E
\]
in which $U/F$ is unramified, $T/U$ has degree prime to $p$ (and its
normal closure is abelian), and the normal closure of $E/T$ has
$p$-power degree. Thus it suffices to treat these three cases separately.
In the first two cases (unramified, or degree prime to $p$),
we have a valuation basis thanks to Proposition~\ref{P:defect basis}, so
there is nothing more to check. In the third case, an exercise in finite group
theory shows that $E/T$ 
can be written as a tower of $\ZZ/p\ZZ$-extensions, to each of which 
we may apply Proposition~\ref{P:degree p}. This yields the desired result.
\end{proof}

\section{$p$-bases and their discrepancies} \label{sec:p-bases}

\begin{defn}
Let $E/F$ be an extension of fields
of characteristic $p>0$ such that $E^p \subseteq F$, so that
in particular $E/F$ is purely inseparable.
A \emph{$p$-basis} of $F$ over $E$ is a
subset $S$ of $E$ such that the products
$\prod_{s \in S} s^{e_s}$, with $e_s \in \{0, \dots, p-1\}$ and all
but finitely many equal to zero, are all distinct and form a basis $T$ of
$E$ as an $F$-vector space. We call $T$
the \emph{associated basis} of $S$; note that $F^* T^*$ forms a group
under multiplication (where $T^* = T \setminus \{0\}$).
A $p$-basis always exists; this is easily seen in case $[F:E^p]$ is finite
(the case we are interested in), but also turns out to be true in general
by Zornication \cite[Section~A.1.3]{eisenbud}.
\end{defn}

For a $p$-basis of an extension of valued fields, 
one has a filtration vaguely reminiscent of the ramification
filtration on the Galois group of a valued field.

\begin{lemma} \label{L:add el disc}
Let $(F,v)$ be a valued field, and let $E$ be a finite extension of $F$ equipped
with an extension of $v$.
For $w,x \in E$, we have
\[
\ndist(wx, F) \geq \min\{\ndist(w,F), \ndist(x,F)\}.
\]
\end{lemma}
\begin{proof}
For any $y,z \in F$ with $v(y) \geq v(w)$, $v(z) \geq v(x)$, we have
$v(yz) \geq v(wx)$ and 
\begin{align*}
v(wx-yz) - v(wx) &\geq \min\{v(w(x-z)), v((w-y)z)\} - v(w) - v(x) \\
&\geq \min\{v(x-z) - v(x), v(w-y) - v(w)\}.
\end{align*}
The supremum of the last expression is $\min\{\ndist(w,F), \ndist(x,F)\}$,
while the supremum of the first expression is at most $\ndist(wx,F)$. This yields
the claim.
\end{proof}

\begin{defn}
Let $(F,v)$ be a valued field of characteristic $p>0$.
Let $E$ be a finite extension of $F$ such that $E^p \subseteq F$,
viewed as a valued field via the unique extension of $v$.
Let $S$ be a $p$-basis of $E$ over $F$ of cardinality $n$, 
with associated basis $T$. For $r \in \Gamma$, put
\[
U_r= \{t \in T: \ndist(t, F) \geq r\};
\]
by Lemma~\ref{L:add el disc},
$F^* U_r^*$ is a subgroup of $F^* T^*$. Moreover,
$U_0 = T$, $U_r = \cap_{s < r} U_s$ for $r \in \RR_{>0}$, 
and $U_r = \{0\}$ for $r$ sufficiently large.
For $i=1, \dots, n$, define the
\emph{$i$-th normalized distance} from $T$ to $F$, denoted
$\ndist_i(T,F)$, to be the supremum of those $r$ for which 
$\# U_r \geq p^{n-i+1}$.
\end{defn}

\begin{prop} \label{P:p-basis disc}
Let $(F,v)$ be a valued field of characteristic $p>0$.
Let $E$ be an extension of $F$ of degree $p^n$ such that $E^p \subseteq F$,
viewed as a
valued field via the unique extension of $v$.
Let $S$ be a $p$-basis of $E/F$ and let $T$ be its associated basis. Then
\[
d(T) \geq (p-1) \sum_{i=1}^n \ndist_i(T,F).
\]
\end{prop}
\begin{proof}
We first define $x_1, \dots, x_n$ as follows. For $i =n, \dots, 1$,
given a choice of $x_{i+1}, \dots, x_n$, choose $x_i \in T$ to maximize
$\ndist(x_i, F)$, subject to the restriction that $x_i, \dots, x_n$
should be linearly independent in $F^*T^*/F^*$. It follows that
\[
\ndist(x_i, F) =  \ndist_i(T,F) \qquad (i=1, \dots, n).
\]

With this definition in hand, we 
prove the inequality $d(T) \geq (p-1)\sum_{i=1}^n \ndist_i(T,F)$.
For any $y_1, \dots, y_n \in F$ with $v(y_i) \geq v(x_i)$, we have
\[
d(T; \prod_{i=1}^n (x_i - y_i)^{p-1}) = \sum_{i=1}^{n} (p-1) (v(y_i-x_i) - v(x_i)).
\]
The supremum of the right side is $(p-1) \sum_{i=1}^{n} \ndist(x_i, F)
= (p-1) \sum_{i=1}^n \ndist_i(T,F)$, while the supremum of the left side
is at most $d(T)$.
\end{proof}

\begin{remark}
We are not sure whether the inequality in Proposition~\ref{P:p-basis disc}
is always achieved, even for $F$ henselian; 
this question can be reformulated as follows.
Put $E_i = F(x_1, \dots, x_i)$
(with $E_0 = F$), and view $\{1, x_i, \dots, x_i^{p-1}\}$ as a
basis for $E_i$ over $E_{i-1}$; the additivity of discrepancies yields
\[
d(T) = \sum_{i=1}^n d(\{1, x_i, \dots, x_i^{p-1}\}).
\]
By Proposition~\ref{P:degree p}, we then have
\[
d(T) = \sum_{i=1}^n (p-1) \ndist(x_i, E_{i-1}).
\]
Consequently, equality holds in Proposition~\ref{P:p-basis disc} if and only
if
\begin{equation} \label{eq:ndist comparison}
\ndist(x_i, E_{i-1}) = \ndist(x_i, F)
\end{equation}
for $i=1, \dots, n$.
\end{remark}

\end{document}